\documentclass{amsart}
\usepackage{latexsym, amsthm, amssymb}

\usepackage[all]{xy}

\usepackage[usenames, dvipsnames]{color}

\usepackage{mathrsfs} 

\theoremstyle{plain}
\newtheorem{Theorem}[equation]{Theorem}
\newtheorem{Lemma}[equation]{Lemma}
\newtheorem{Corollary}[equation]{Corollary}

\newtheorem*{Conjecture*}{Conjecture}
\newtheorem{Proposition}[equation]{Proposition}

\theoremstyle{definition}

\newtheorem{Paragraph}[equation]{}

\theoremstyle{remark}
\newtheorem{Remark}[equation]{Remark}

\numberwithin{equation}{section}

\newcommand{\h}{{\text{\rm h}}}

\newcommand{\F}{{\mathbb F}}

\renewcommand{\H}{{\text{\rm H}}}

\renewcommand{\O}{{\text{\rm O}}}

\newcommand{\Q}{{\mathbb Q}}
\newcommand{\R}{{\mathbb R}}
\renewcommand{\S}{{\mathcal S}}

\newcommand{\Z}{{\mathbb Z}}

\newcommand{\car}{{\text{\rm char}}}

\newcommand{\cok}{{\text{\rm cok}}}

\newcommand{\cont}{{\text{\rm cont}}}
\newcommand{\cs}{{\text{\rm cs}}}

\newcommand{\df}{{\,\overset{\text{\rm df}}{=}\,}}

\renewcommand{\div}{{\text{\rm div}}}
\renewcommand{\dim}{{\text{\rm dim}}}

\newcommand{\et}{{\text{\rm \'et}}}

\newcommand{\im}{{\text{\rm im}}}
\newcommand{\ind}{{\text{\rm ind}}}
\renewcommand{\inf}{{\text{\rm inf}\,}}

\newcommand{\isim}{{\;\overset{\sim}{\longrightarrow}\;}}
\newcommand{\isom}{{\;\simeq\;}}

\newcommand{\lcm}{{\text{\rm lcm}}}

\newcommand{\pf}{{\noindent{\it Proof.}\;\;}}

\newcommand{\red}{{\text{\rm red}}}

\newcommand{\res}{{\text{\rm res}}}
\newcommand{\rk}{{\text{\rm rk}}}

\newcommand{\sh}{{\text{\rm sh}}}

\newcommand{\Ass}{{\text{\rm Ass}}}

\newcommand{\Br}{{\text{\rm Br}}}

\newcommand{\Div}{{\text{\rm Div}\,}}

\newcommand{\Frac}{{\text{\rm Frac}\,}}

\newcommand{\Hom}{{\text{\rm Hom}}}

\newcommand{\Pf}{{\noindent{\it Proof.}\;\;}}

\newcommand{\Supp}{{\text{\rm Supp}}}
\newcommand{\Spec}{{\text{\rm Spec}\,}}

\begin{document}

\title[Tame Covers and Cohomology of Relative Curves]
{Tame Covers and Cohomology of Relative Curves
over Complete Discrete Valuation Rings
with Applications to the Brauer Group
}

\author{Eric Brussel and Eduardo Tengan}
\address 
{Department of Mathematics \& Computer Science\\
Emory University\\
Atlanta, GA 30322\\ USA}
\email{brussel@mathcs.emory.edu}
\address
{Instituto de Ci\^encias Matem\'aticas e de Computa\c c\~ao\\
Universidade de S\~ao Paulo\\
S\~ao Carlos, S\~ao Paulo\\ Brazil}
\email
{etengan@icmc.usp.br}

\subjclass{11G20, 11R58, 12G05, 16K20, 16K50}




\begin{abstract}
We prove the existence of noncrossed product division algebras and indecomposable
division algebras of unequal period and index over the function field of any $p$-adic
curve, extending the results and methods of \cite{BMT}.
\end{abstract}

\maketitle

\section{Introduction}

We study the cohomology and the Brauer group of a field $F$ that is finitely generated
and of transcendence degree one over the $p$-adic field $\Q_p$.
Such a field is always the function field of a regular (projective, flat) relative curve $X/\Z_p$.
In \cite{BMT} it was shown that if $F$ admits a {\it smooth} model
$X/\Z_p$ then there exist noncrossed product $F$-division algebras,
and indecomposable $F$-division algebras of unequal prime-power period and index.
These were constructed from objects defined over the generic point $(p)$ of the closed
fiber $X_0=X\otimes_{\Z_p}\F_p$ using a
homomorphism $\lambda:\Br(F_p)'\to\Br(F)$ that splits the restriction map
$\res:\Br(F)\to\Br(F_p)$, where the subscript $p$ denotes completion with respect to 
the valuation defined by $(p)$, and ``$\;'\;$'' denotes ``prime-to-$p$ part''. 
The fields $F$ that have a smooth $X/\Z_p$ include fields such as $\Q_p(t)$ 
but do not include the function
fields of {\it all} $p$-adic curves.

In this paper we generalize the machinery and results of \cite{BMT} to arbitrary $p$-adic curves.
We prove that if $F$ is the function field of a $p$-adic curve then there exist
noncrossed product $F$-division algebras, and indecomposable $F$-division algebras
of unequal prime-power period and index.  The machinery we develop here is used in \cite{BT11b}
to prove that every $F$-division algebra of prime period $\ell$ and index $\ell^2$ decomposes 
into two cyclic $F$-tensor factors, 
hence is a crossed product, generalizing Suresh's result \cite{Sur10},
which assumes roots of unity.  
In the terminology of \cite[Sections 3,4]{ABGV},
this shows the $\Z/\ell$-length in $\H^2(F,\mu_\ell)$ equals the $\ell$-Brauer dimension,
which is two by a theorem of Saltman (\cite[Theorem 3.4]{Sa97}).
In general our work is motivated by the work of Saltman over these fields in \cite{Sa97} and \cite{Sa07}
(see also \cite{Br10}).

We summarize the technical results.
Let $R$ be a complete discrete valuation ring with fraction field $K$,
and let $F$ be a finitely generated field extension of $K$ of transcendence degree one.
Let $X/R$ be a regular (projective, flat) model for $F$
whose closed fiber $X_0$ has normal crossings on $X$.
Let $C=X_{0,\red}$, let $\{C_i\}$ denote the set of irreducible components of $C$,
let $\S$ be the set of singular points of $C$, 
and set $F_C=\prod_i F_{C_i}$, the product of
the completions with respect to the valuations on $F$ defined by the $C_i$.
We construct for any prime-to-$\car(\kappa(C))$ number $n$, 
any integer $r$, and any $q\geq 0$, a homomorphism
$$
\lambda:\H^q(\O_{C,\S},\Z/n(r))\to\H^q(\O_{X,\S},\Z/n(r))\to\H^q(F,\Z/n(r))
$$
that splits the restriction map $\H^q(\O_{X,\S},\Z/n(r))\to\H^q(\O_{C,\S},\Z/n(r))$.
We use $\lambda$ to construct another map $\lambda$ that splits the subgroup
of the image of $\H^q(F,\Z/n(r))\to\H^q(F_C,\Z/n(r))$
consisting of tuples of classes that are unramified at $\S$ and glue along $\S$.
When $q=2$, $\Z/n(r)=\mu_n$, and $R=\Z_p$, we show that $\lambda$ preserves index.
This allows us to construct indecomposable
$F$-division algebras and noncrossed product $F$-division algebras as mentioned above,
in the same manner as \cite{BMT}.
When the dual graph of the closed fiber $X_0$ has nontrivial topology, i.e., 
nonzero (first) Betti number, we construct cyclic covers of $X$
that are (defined and) trivial at every point of $X$ except the generic point of $X$.
These arise as cyclic covers of the closed fiber $X_0$ that are trivial at every point, and transported
to $X$ via $\lambda$.
When $R=\Z_p$ they are the completely split cyclic covers of 
Saito (\cite{Sai85}).
We thank Colliot-Th\'el\`ene for drawing our attention to these interesting specimens.

\section{Background and Conventions}

We use \cite[Chapter 8,9]{Liu},
\cite[Section 2]{GM}, and \cite[Chapter XIII]{SGA1} for many of the following definitions.

\Paragraph{General Conventions.}
Let $S$ be an excellent scheme, $n$ a number that is invertible on $S$,
and $\Lambda=(\Z/n)(i)$ the twisted \'etale sheaf.
We write $\H^q(S,\Lambda)$ for the \'etale cohomology group,
and if $\Lambda$ is understood (or doesn't matter)
we write $\H^q(S,r)$ instead of $\H^q(S,\Lambda(r))$,
and $\H^q(S)$ in place of $\H^q(S,0)$.
If $S=\Spec A$ for a ring $A$ then we 
write $\H^q(A,r)$.
If $T$ is an integral scheme contained in $S$ we write $\kappa(T)$ for its function field. 
If $T\to S$ is a morphism of schemes then
the restriction $\res_{T|S}:\H^q(T)\to\H^q(S)$ is defined,
and we write $\beta_S=\res_{T|S}(\beta)$.  
If $Z\subset S$ is a subscheme, we write $Z_T$ for the preimage $Z\times_S T$.

\Paragraph{Valuation Theory.}
If $v$ is a valuation on a field $F$ we write $\kappa(v)$ for
the residue field of the valuation ring $\O_v$, and $F_v$ for the completion of $F$
at $v$.  
If $S$ is a connected normal scheme with function field $F$ 
and $v$ arises from a prime divisor $D$ on $S$,
we write $v_D$ for $v$, $\kappa(D)$ for $\kappa(v)$, and $F_D$ for $F_v$.
If $D$ is a sum of prime divisors $D_i$ we write
$F_D=\prod_i F_{D_i}$.
Each $f\in F^*$ defines a divisor $\div(f)=\sum v_D(f)D$,
where the (finite) sum is over prime divisors on $S$.

Recall that if $F=(F,v)$ is a discretely valued field and $\alpha\in\H^q(F,\Lambda)$
then $\alpha$ has a {\it residue} $\partial_v(\alpha)$ in $\H^{q-1}(\kappa(v),\Lambda(-1))$.
%
More generally,
suppose $T$ is a noetherian scheme,
$\xi$ is a generic point of $T$, and $\alpha\in\H^q(T,\Lambda)$.
Then for each discrete valuation $v$ on the field $F=\kappa(\xi_\red)$
$\alpha$ has a residue 
$$
\partial_v(\alpha)\df\partial_v(\alpha_F)\in\H^{q-1}(\kappa(v),\Lambda(-1))
$$
We say $\alpha$ is {\it unramified at $v$} if $\partial_v(\alpha)=0$, {\it ramified at $v$}
if $\partial_v(\alpha)\neq 0$, and {\it tamely ramified at $v$} if $\partial_v(\alpha)$
is contained in the prime-to-$\text{\rm char}(\kappa(v))$ part of $\H^{q-1}(\kappa(v),\Lambda(-1))$.
If $\alpha$ is unramified at $v$ the {\it value} of $\alpha$ at $v$ is the
element 
$$
\alpha(v)=\res_{F|F_v}(\alpha_F)\in\H^q(\kappa(v),\Lambda)
\leq\H^q(F_v,\Lambda)
$$
(see \cite[7.13, p.19]{GMS}).
Suppose $T\to S$ is a birational morphism of noetherian schemes (see \cite[I.2.2.9]{EGAI}).
The {\it ramification locus of $\alpha$ on $S_\red$}
is the sum of the prime divisors on $S_\red$ that determine valuations at
which $\alpha$ is ramified, over all generic points of $S_\red$.

Let $D$ be a divisor on a noetherian normal
scheme $S$, set $U=S-D$, and for each generic point $\xi$ of $\Supp\,D$, let $K_\xi$
denote the fraction field of the discrete valuation ring $\O_{S,\xi}$.
We say a morphism $\rho:T\to S$
is {\it tamely ramified along $D$} if $T$ is normal, $\rho_U:V=T\times_S U\to U$ is \'etale,
and for each generic point $\xi$ of $\Supp\,D$, the \'etale $K_\xi$-algebra $L$ defined by
$\Spec L=V\times_U\Spec K_\xi$ is tamely ramified with respect to $\O_{S,\xi}$.
Since $L/K_\xi$ is \'etale it is a finite product of separable field extensions
of $K_\xi$, and $L$ is tamely ramified if each field extension is tamely ramified 
(with respect to $\O_{S,\xi}$) in the usual sense.
If $S$ is a noetherian scheme whose irreducible components are normal, we'll
say a morphism $\rho:T\to S$ is tamely ramified along $D$ if again $V\to U=S-D$ is \'etale,
and the restriction $\rho_{S_i}$
to each irreducible component $S_i$ of $S$ is tamely ramified along $D_{S_i}$.
If $S=\Spec A$ and $T=\Spec B$, we will also say $B$ is a {\it tamely ramified $A$-algebra}.
We say a map $\rho:T\to S$ of noetherian schemes 
is a {\it cover} if it is finite, generically \'etale, 
and each connected component of $T$ dominates a connected component of $S$.



\Paragraph{Relative Curves.}\label{setup}
In this paper, $R$ will be a complete discrete valuation ring with residue field $k$
and fraction field $K$, 
$F$ will be a field finitely generated of transcendence degree one over $K$,
and $X/R$ will be
a regular 2-dimensional scheme $X$ that is flat and projective over $\Spec R$
and has function field $K(X)=F$.  We call $X/R$ a {\it regular relative curve}, 
write $X_0=X\otimes_R k$ for the closed fiber,
$C=X_{0,\red}$ for the reduced subscheme underlying the closed fiber, 
and $C_1,\dots,C_m$ for the irreducible components of $C$.
We assume each $C_i$ is regular, and at most two of them meet at any closed point of $X$,
a situation that can always  be obtained by blowing up using Lipman's embedded resolution theorem
(see \cite[9.2.4]{Liu}).
For all closed points $z\in X$, we have $\dim\O_{X,z}=2$ by \cite[8.3.4(c)]{Liu},
and since $X$ is regular, $\O_{X,z}$ is factorial by Auslander-Buchsbaum's theorem.

We say an effective divisor $D$ on a relative curve $X/R$ is {\it horizontal}
if each of its irreducible components maps surjectively (hence finitely) to $\Spec R$,
and {\it vertical} if its support is contained in the support of the closed fiber.
If $D$ is a reduced and irreducible horizontal divisor then it is flat over 
$\Spec R$, since $R$ is a discrete valuation ring.
Every effective divisor on a 
relative curve $X/R$ is a sum of horizontal and vertical divisors,
and the horizontal prime divisors are exactly the closures of the closed points of the generic fiber
(\cite[8.3.4(b)]{Liu}).
Since $R$ is henselian, each irreducible horizontal divisor has a single closed point.

\Paragraph{Distinguished Divisors.}\label{distinguished}
In general there will be many horizontal divisors on a relative curve $X$
that restrict to a given divisor on $C$.  
In order to construct our lifts of covers and cohomology classes
from $C$ to $X$ we select a single regular horizontal divisor for each closed point, as follows.

\Proposition\label{distingprop}
Assume the setup of \eqref{setup}.
Then for each closed $z\in X$ there exists a regular irreducible horizontal divisor $D\subset X$
that intersects each irreducible component of $C$ passing through $z$ transversally at $z$.
\rm

\begin{proof}
Transversal intersection with a single component is by
\cite[8.3.35(g)]{Liu} and its proof (see also \cite[21.9.12]{EGAIV:d}).
Thus if $z\in C_i\cap C_j$ ($i\neq j$) and $t_i$ and $t_j$ are local equations for $C_i$ and $C_j$,
then we have local equations $f_i$ and $f_j$ for effective regular horizontal divisors
such that $(f_i,t_i)=(f_j,t_j)=\frak m_z\subset\O_{X,z}$.
If $(f_j,t_i)=\frak m_z$ or $(f_i,t_j)=\frak m_z$ then a suitable $D$ is defined locally
by $f_j$  or $f_i$.  Otherwise $(f_i+f_j,t_i)=(f_i+f_j,t_j)=\frak m_z$,
and we define $D$ locally by $f_i+f_j$.
The rest of the proof proceeds as in \cite[3.3.35]{Liu}. 
\end{proof}\rm

We fix a set of these (prime) divisors, 
and let $\mathscr D$ denote the set of supports of the semigroup they generate in $\Div X$.
We will say a divisor $D$ is {\it distinguished}
and write $D\in\mathscr D$ whenever $D$ is reduced and supported in $\mathscr D$.
Though $\mathscr D$ is fixed in principle, we reserve the right to declare any
divisor satisfying the definition to be a member of this set retroactively.
Let $\mathscr D_\S$ denote the subset that {\it avoids $\S$}.
Note that each $D\in\mathscr D$ is a disjoint union of its irreducible components,
each of which meets each irreducible component of $C$
transversally.



\section{Structure of Tame Covers}

\Lemma[Structure]\label{structure}
Assume the setup of (\ref{setup}).
Suppose $\rho:Y\to(X,D)$ is a tamely ramified cover, where $D\in\mathscr D$.
Then 
\begin{enumerate}
\item[a)]
The structure map $\rho:Y\to X$ is flat.
\item[b)]
$Y/R$ is a regular relative curve, $Y_{0,\red}=C_Y$,
each irreducible component of $C_Y$ is regular,
$\S_Y$ is the set of singular closed points of $C_Y$,
and exactly two irreducible components of $C_Y$ meet at each point of $\S_Y$.
\item[c)]
The support of the irreducible components of $D'_Y$ for $D'\in\mathscr D$
generate a set $\mathscr D_Y$ of distinguished divisors on $Y$.
\end{enumerate}
\rm

\pf
Since $Y\to X$ is finite, $\dim(X)=\dim(Y)=2$ by \cite[5.4.2]{EGAIV:b}, and $Y\to\Spec R$
is projective as the composition of projective morphisms (\cite[3.3.32]{Liu}).
Let $y\in Y$ be a closed point and set $x=f(y)$, $A=\O_{X,x}$, $B'=\O_{Y,x}$, and $B=\O_{Y,y}$.
Choose a geometric point over $x$ that lifts to each point of $Y$ lying over $x$,
and use this in the following to define the strict henselizations with respect to the maximal ideals
of these points.

Since the statements involving $D$ are local and $D$ is a disjoint union of its irreducible
components we may assume $D$ is irreducible.
Let $C_i\subset C$ be a (regular) irreducible component going through $x$,
and let $\{f,t\}\subset A$ be the regular system of parameters formed by local equations for 
the distinguished prime divisor passing through $x$, and for $C_i$, respectively.
Then the strict henselization $A^\sh$ of $A$ with respect to the maximal ideal of $A$
is a two-dimensional regular local ring, faithfully flat over $A$, with regular system
of parameters $\{f,t\}$ (see \cite[18.8]{EGAIV:d}).

If $x\not\in D$ then $B'\otimes_A A^\sh$ is a finite \'etale $A^\sh$-algebra
by base change, since $\rho|_{X-D}$ is finite-\'etale.
If $x\in D$ then $B'\otimes_A A^\sh$ is a finite tamely ramified $A^\sh$-algebra
by \cite[Lemma 2.2.8]{GM}.
By \cite[18.8.10]{EGAIV:d}, $B^\sh$ is a factor of the direct product decomposition of
$B'\otimes_A A^\sh$, hence 
$B^\sh$ is a finite tamely ramified local $A^\sh$-algebra, 
in particular it is a normal local ring, hence it is a normal domain.
It follows that
$B^\sh$ is the integral closure of $A^\sh$ in the field $\widetilde L\df\Frac B^\sh$.
Since the tame fundamental group of the strictly henselian regular local ring $A^\sh$ 
is abelian (\cite[XIII.5.3]{SGA1})
the field extension $\widetilde L/\Frac A^\sh$ is Galois,
and by Abhyankar's Lemma (\cite[A.I.11]{FK}, see also \cite[Corollary 2.3.4]{GM})
$$
B^\sh=A^\sh[T]/(T^e-f)\quad(\text{some }e\geq 1)
$$ 
By \cite[Lemma 1.8.6]{GM}
$B^\sh$ is a regular (2-dimensional) local ring with system of parameters $\{\root e\of f,t\}$.
Since $B\to B^\sh$ is faithfully flat
and $B^\sh$ is regular, $B$ is regular by flat descent (\cite[6.5.1]{EGAIV:b} or \cite[23.7(i)]{Mat}), 
and since $B$ is the local ring of an arbitrary closed point, we conclude $Y$ is regular.
It follows that $\rho:Y\to X$ is flat by \cite[23.1]{Mat}, proving (a), and since $Y$ is regular
and $Y\to\Spec R$ is flat and projective, $Y/R$ is a regular relative curve.

We derive a system of parameters for $B$.
The prime ideal $(\root e\of f)\subset B^\sh$ is the only one lying over $(f)A^\sh$ since,
for $\kappa(f)=\Frac A^\sh/(f)A^\sh$,
the ring $B^\sh\otimes_{A^\sh}\kappa(f)=\kappa(f)[T]/(T^e)$ of the fiber over $\Spec\kappa(f)$ 
consists of a single prime ideal.
The image $(\root e\of f)$ in $\Spec B$ is therefore a unique prime $(g)\subset B$ lying
over $(f)\subset A$, and $(\root e\of f)$ is the unique prime lying over $(g)$.
Therefore, since $B\to B^\sh$ is unramified, $(g)B^\sh=(\root e\of f)$.
Since $B\to B^\sh$ is faithfully flat, $IB^\sh\cap B=I$ for all ideals $I$ of $B$ (by e.g.
\cite[Exercise 3.16]{AM}), so since $(g,t)B^\sh=(\root e\of f,t)$ is maximal, 
$(g,t)B^\sh\cap B=(g,t)$ is the maximal ideal of $B$.
Thus $\{g,t\}$ is a regular system for $B$.

Since $t$ is a local equation for $\rho^{-1}C_i$, $\rho^{-1}C_i$ is regular
and irreducible at $y$ for each $C_i$ passing through $x$.
In particular $C_Y=\bigcup_i\rho^{-1} C_i$ is reduced, and so equals $Y_{0,\red}$.  
Since at most two irreducible components of $C$ meet at $x$, the same holds for $C_Y$ at $y$,
and $y$ is a singular point on $C_Y$ if and only if $x=f(y)\in\S$.
This completes the proof of (b).

If $D'\in\mathscr D$ is the distinguished (horizontal) prime divisor running through $x$
then there is a single irreducible component of $D_Y'$ passing through $y$, whose 
support $D_{Y,\red}'$ has local equation $g$ at $y$.
Thus each irreducible component of $D_Y'$ covers $D'$, hence $\Spec R$,
hence $D_Y'$ is horizontal.
Since $g$ is part of the regular system $\{g,t\}$ at the arbitrary closed point $y$
we see that $D_{Y,\red}'$ is regular, and since $t$
is a local equation for an arbitrary irreducible component of $C_Y$ passing through $y$,
$D_{Y,\red}'$ intersects each component of $C_Y$ transversally.
Thus the support of the irreducible components of $D_Y'$ generate a set
of distinguished divisors $\mathscr D_Y$ for $Y$.  This proves (c).



\hfill $\blacksquare$

\Lemma\label{covers}
Suppose $X$ is a regular noetherian scheme and
$L$ is an \'etale $K(X)$-algebra that is tamely ramified along a divisor $D$.
Then the normalization $Y$ of $X$ in $L$ defines a tamely ramified cover
$\rho:Y\to(X,D)$.
\rm

\Pf
Since $X$ is regular, its connected components are integral regular schemes, 
hence we may assume $X$ is integral.
Since $L/K(X)$ is \'etale, $L$ is a product of finite separable field extensions of $K(X)$,
hence we may assume $L/K(X)$ is itself a finite separable field extension.
Then the normalization $Y$ exists, $Y$ is normal by definition,
and $\rho:Y\to X$ is finite by \cite[4.1.25]{Liu}.
Since $Y$ is normal and connected it is irreducible, so $Y$ dominates $X$.
Let $U=X-D$, and set $V=Y\times_X U$.
Since $X$ is normal, $Y$ is connected, and $\rho|_V$ is unramified, 
$\rho|_V$ is \'etale by \cite[I.9.11]{SGA1} (see also \cite[I.3.20]{M}).
Therefore $Y\to (X,D)$ is a tamely ramified cover.

\hfill $\blacksquare$

The next lemma shows how distinguished divisors split in tamely ramified covers.

\Lemma\label{lemma2}
Assume the setup of (\ref{setup}).
Suppose $\rho:Y\to (X,D)$ is a tamely ramified cover, where $D\in\mathscr D$,
and $D'\in\mathscr D_\S$ is irreducible.
Suppose $E\subset D_{Y,\red}'$ is a distinguished prime divisor lying over $D'$ as in
Lemma~\ref{structure}(c),
$y= E\times_Y C_Y$, and $x=D'\times_X C$.
Then $y$ and $x$ are regular closed points,
and the ramification (resp. inertia) degree of $v_E$ over $v_{D'}$ 
equals the ramification (resp. inertia)
degree of $v_y$ over $v_x$.
\rm

\Pf
Since we assume \eqref{setup} and $D\in\mathscr D$ we have Lemma~\ref{structure},
which shows $C_Y$ is reduced and $E\subset D_{Y,\red}'$ is distinguished.
Note that either $D'\cap D=\varnothing$ or $D'\subset D$.
Since $D'$ and $E$ are distinguished and avoid the singular points of $X$ and $Y$,
they intersect the reduced closed fibers $C$ and $C_Y$ transversally, hence
$x=D'\times_X C$ and $y=E\times_Y C_Y$ are regular closed points.
We must show that $[\kappa(E):\kappa(D')]=[\kappa(y):\kappa(x)]$
and that $v_E(f)=v_y(f_0)$, where $f\in\O_{X,D'}$ is a local equation for $D'$ on $X$
and $f_0\in\O_{C,x}$ is a local equation for $x$ on $C$. 
 
Since $D'$ is horizontal and irreducible, 
$D'=\Spec S$ for $S$ a finite local $R$-algebra by \cite[I.4.2]{M},
and $S$ is a discrete valuation ring since $D'$ is regular.
The map $E\to\rho^{-1}D'\to D'$ is finite as a composition of finite morphisms, hence
$E=\Spec T$ for $T$ a finite local $S$-algebra, again a discrete valuation ring since
$E$ is regular.
Since $S$ is a discrete valuation ring and $S\to T$ is finite, $T$ is a free $S$-module of finite
rank, and so $[T:S]$ is well defined.
Since the generic point of $E$ lies over that of $D'$, we have $\Frac T=T\otimes_S\Frac S$,
hence $[\kappa(E):\kappa(D')]=[T:S]$.

Let $A=\O_{X,x}$, $B=\O_{Y,y}$,
let $t$ be a local equation for $C$ at $x$, and
set $A_0=A/(t)$ and $B_0=B/(t)B$, the (reduced) local rings of the fibers 
through $x$ and $y$, as in the proof of Lemma~\ref{structure}.
Already $\kappa(x)=S\otimes_A A_0$ and $\kappa(y)=T\otimes_B B_0$
by the transversality of the intersections.
Since $B_0=B\otimes_A A_0$ we have
$\kappa(y)=T\otimes_A A_0$,
hence $[\kappa(y):\kappa(x)]=[T:S]=[\kappa(E):\kappa(D')]$ by base change.

Let $g\in A$ be defined as above.
To compute the ramification degree, note that since $B\to B^\sh$ is faithfully flat, 
$(g^e)B=(g^e)B^\sh\cap B=(f)B^\sh\cap B=(f)B$,
hence $g^e=fu$ for some $u\in B^*$.
Since $f$ and $g$ are uniformizers for $v_{D'}$ and $v_{E}$, respectively, it follows that
$e(v_{E}/v_{D'})=v_{E}(f)=e$.
On the other hand, let $f_0$ be the image of $f$ in $A_0$, and let
$g_0$ be the image of $g$ in $B_0$. 
Then $f_0$ cuts out the closed point $x$ on $C$
and $g_0$ cuts out $y$ on $C_Y$ by transversality.
Thus $f_0$ and $g_0$ are uniformizers for $v_x$ and $v_y$, and since $g_0^e=f_0 u_0$,
where $u_0$ is the image of $u$ in $B_0^*$,
we have $e(v_y/v_x)=v_y(f_0)=e$, as desired.
This completes the proof.

\hfill $\blacksquare$

\section{Lifting Cohomology Classes}

\Paragraph\label{setup2}
Let $k$ be a field, and
let $C/k$ be a reduced connected projective curve with regular irreducible components $C_1,\dots,C_m$,
at most two of which meet at any closed point.
Denote the singular points of $C$ by $\S$
and write $\O_{C,\S}$ for the semilocal ring $\varinjlim_U\O_C(U)$,
where $U$ varies over (dense) open subsets of $C$ containing $\S$.
Then $\O_{C,\S}$ is a subring of the rational function ring $\kappa(C)=\prod_i \kappa(C_i)$.
For each $z\in \S\cap C_i$, let $K_{i,z}=\Frac\O_{C_i,z}^\h$, a field since $z$ is a normal point of
$C_i$, and if $\alpha_i\in\H^q(\kappa(C_i))$,
let $\alpha_{i,z}$ denote the image in $\H^q(K_{i,z})$.

\Lemma[Gluing]\label{gluing}
Assume the setup of (\ref{setup2}). 
There exists an element $\alpha\in\H^q(\O_{C,\S},\Lambda)$
that restricts to
$\alpha_C=(\alpha_1,\dots,\alpha_m)\in\bigoplus_i\H^q(\kappa(C_i),\Lambda)$
if and only if $\alpha_i$ is unramified at each $z\in\S\cap C_i$,
and $\alpha_{i,z}=\alpha_{j,z}$ whenever $z\in C_i\cap C_j$.
\rm

\Pf
There is an exact sequence (\cite[III.1.25]{M})
\begin{align*}
0
\longrightarrow \H^0_\S(\O_{C,\S})&\longrightarrow\H^0(\O_{C,\S})\longrightarrow
\H^0(\kappa(C))\longrightarrow\H^1_\S(\O_{C,\S})\longrightarrow\\
&\quad\longrightarrow\H^1(\O_{C,\S})\longrightarrow\H^1(\kappa(C))
\longrightarrow\H^2_\S(\O_{C,\S})\longrightarrow\\
&\qquad\longrightarrow\H^2(\O_{C,\S})\longrightarrow\H^2(\kappa(C))
\longrightarrow\H^3_\S(\O_{C,\S})
\tag{$*$}
\end{align*}
where the maps into the direct sum are restrictions.
Since $\S$ is a disjoint union of closed points,
$\H^q_\S(\O_{C,\S})=\bigoplus_{z\in\S}\H^q_z(\O_{C,\S})=\bigoplus_{z\in \S}\H^q_z(\O_{C,z}^\h)$ 
by excision (\cite[III.1.28, p.93]{M}).
Since $\Lambda$ is a smooth group scheme, $\H^q(\O_{C,z}^\h)=\H^q(\kappa(z))$,
by the cohomological Hensel's lemma \cite[III.3.11(a), p.116]{M}.
Since the $C_k$ are regular and at most two of them meet at any $z\in\S$, 
we have
$\Spec\O_{C,z}^\h-\{z\}=\Spec(K_{i,z}\times K_{j,z})$ for some $i$ and $j$,
and an ``excised'' exact sequence
\begin{align*}
0
\longrightarrow \H^0_z(\O_{C,\S})&\longrightarrow\H^0(\kappa(z))
\longrightarrow\H^0(K_{i,z}\times K_{j,z})\longrightarrow\H^1_z(\O_{C,\S})\longrightarrow\\
&\quad\longrightarrow\H^1(\kappa(z))\longrightarrow 
\H^1(K_{i,z}\times K_{j,z})\longrightarrow\H^2_z(\O_{C,\S})\longrightarrow\\
&\qquad \longrightarrow\H^2(\kappa(z))\longrightarrow \H^2(K_{i,z}\times K_{j,z})\longrightarrow\H^3_z(\O_{C,\S})\longrightarrow\cdots\\
\end{align*}
where the map $\H^q(\kappa(z))\to\H^q(K_{i,z}\times K_{j,z})=
\H^q(K_{i,z})\oplus\H^q(K_{j,z})$ 
is the diagonal map given by inflation from $\kappa(z)$ to the ``local fields'' $K_{i,z}$ and $K_{j,z}$.
Since $n$ is prime-to-$p$, the map $\H^0(\kappa(z))\to\H^0(K_{i,z})$ is an isomorphism,
so $\H^0_z(\O_{C,\S})=0$, and
for $q\geq 1$ we have short exact Witt-type sequences
\[
0\to\H^q(\kappa(z))\to\H^q(K_{i,z})\xrightarrow{\;\partial_z\,}\H^{q-1}(\kappa(z),-1)\to 0
\]
Thus the long exact sequence breaks up into short exact sequences 
\begin{equation}\label{ses}
0\to\H^q(\kappa(z))\to\H^q(K_{i,z}\times K_{j,z})\to\H^{q+1}_z(\O_{C,\S})\to 0
\qquad (q\geq 0)
\end{equation}
By the compatibility of the localization sequence with the excised sequence, the map
$\H^q(\kappa(C_i))\to\H^{q+1}_z(\O_{C,\S})\leq\H^{q+1}_\S(\O_{C,\S})$ of $(*)$ factors through
$\res_{\kappa(C_i)|K_{i,z}}$.
Therefore an element $\alpha_C=(\alpha_1,\dots,\alpha_m)\in\H^q(\kappa(C))$
maps to zero in $\H_\S^{q+1}(\O_{C,\S})$
if and only if each couple $(\alpha_{i,z},\alpha_{j,z})$ is in the image of some $\bar\alpha\in\H^q(\kappa(z))$;
i.e., $\alpha_{i,z}=\alpha_{j,z}$, and both are unramified.
Thus by the exactness of $(*)$, 
$\alpha_C$ comes from $\H^q(\O_{C,\S})$ if and only if each $\alpha_i$ is unramified at each $z\in\S\cap C_i$,
and $\alpha_{i,z}=\alpha_{j,z}$ whenever $z\in C_i\cap C_j$.

\hfill $\blacksquare$

Suppose $C$ is as in \eqref{setup2}.
Since exactly two irreducible components meet at any $z\in\S$
the {\it dual graph} $G_C$ is defined, and consists of a vertex for each irreducible
component of $C$ and an edge for each singular point, 
such that an edge and a vertex are incident when the corresponding singular point
lies on the corresponding irreducible component (\cite[2.23]{Sai85}, see also \cite[10.1.48]{Liu}).
The (first) Betti number for $G_C$ is $\beta_C\,\df\,\rk(\H_1(G_C,\Z))=N+E-V$,  
where $V,E$ and $N$ are the numbers of vertices, edges, and connected components of $G_C$, 
respectively.

\Lemma\label{injects}
Assume the setup of (\ref{setup2}).
Then:
\begin{enumerate}
\item[a)]
For any integer $r$,
$\H^1(C,\Z/n(r))\to\H^1(\O_{C,\S},\Z/n(r))$ is injective.
\item[b)]
The map $\H^q(\O_{C,\S},\Z/n(q-1))\to\H^q(\kappa(C),\Z/n(q-1))$ is injective
for $q=0,2$, and for $q=1$ we have
$$\H^1(\O_{C,\S},\Z/n)\;\isom\;(\Z/n)^{\beta_C}\oplus \Gamma$$
where 
$(\Z/n)^{\beta_C}$ is the kernel of $\H^1(\O_{C,\S},\Z/n)\to\H^1(\kappa(C),\Z/n)$,
and $\Gamma\leq \H^1(\kappa(C),\Z/n)$ 
is the group of tuples that glue as in Lemma~\ref{gluing}.
\end{enumerate}
\rm

\Pf
We suppress the notation for $\Lambda=\Z/n(r)$.
Let $z\in C-\S$ be a (regular) closed point, and set $U=C-\{z\}$, 
a dense open subset containing $\S$.  
The localization exact sequence is
\begin{align*}
0\longrightarrow\H^0_z(C)\longrightarrow&\H^0(C)\longrightarrow\H^0(U)
\longrightarrow\cdots\\
&\cdots\longrightarrow\H^q_z(C)\longrightarrow\H^q(C)\longrightarrow\H^q(U)\longrightarrow\H^{q+1}_z(C)\longrightarrow\cdots
\end{align*}
By excision we have an exact sequence
\begin{align*}
0\longrightarrow\H^0_z(C)\longrightarrow&\H^0(\O_{C,z}^\h)
\longrightarrow\H^0(K_z)\longrightarrow\cdots\\
&\cdots\longrightarrow\H^q_z(C)
\longrightarrow\H^q(\O_{C,z}^\h)\longrightarrow\H^q(K_z)
\longrightarrow\H^{q+1}_z(C)\longrightarrow\cdots
\end{align*}
where $K_z=\Frac\O_{C,z}^\h$.
Since $z$ is a regular point $\O_{C,z}^\h$ is a discrete valuation ring,
and by \cite[Section 3.6]{CT95}
we may replace $\H^{q+1}_z(C)$ with $\H^{q-1}(\kappa(z),-1)$,
and the map from $\H^q(U)$, which factors
through $\H^q(K_z)$, is then the residue map $\partial_z$.
We conclude 
$\H^0(C)=\H^0(U)$, and we have a long exact sequence
\begin{align*}
0\longrightarrow&\H^1(C)\longrightarrow\H^1(U)\xrightarrow{\;\partial_z\;}
\H^0(\kappa(z),-1)\longrightarrow \cdots\\
&\cdots\longrightarrow\H^q(C)\longrightarrow\H^q(U)
\xrightarrow{\;\partial_z\;}\H^{q-1}(\kappa(z),-1)\longrightarrow\cdots
\tag{$**$}
\end{align*}
As $\H^1(\O_{C,\S})=\varinjlim_U\H^1(U)$,
where the limit is over dense open subsets of $C$ containing $\S$,
$\H^1(C)\to\H^1(\O_{C,\S})$ is injective by
the exactness of the injective limit functor, proving (a).

For (b) we go back to $\Lambda=\Z/n$.    
By ($*$) we have an exact sequence
\begin{align*}
0\longrightarrow\H^0(\O_{C,\S},\Z/n)\xrightarrow{\,\phi_1\,}\,
&\H^0(\kappa(C),\,\Z/n)\xrightarrow{\,\phi_2\,}\,\H^1_\S(\O_{C,\S}, \Z/n)\xrightarrow{\,\phi_3\,}\,\\
&\xrightarrow{\,\phi_3\,}\,\H^1(\O_{C,\S},\Z/n)\xrightarrow{\,\phi_4\,}\,\H^1(\kappa(C),\Z/n)
\end{align*}
The groups $\H^0(\O_{C,\S},\Z/n)$ and $\H^0(\kappa(C),\Z/n)$
are finite free $\Z/n$-modules whose ranks are the number of $C$'s connected components $N$
and irreducible components $m$, respectively.
We claim $\H^1_\S(\O_{C,\S},\Z/n)$ is a finite free $\Z/n$-module.
For by (\ref{ses}), for each $z\in\S$ we have an exact sequence
$$
0\longrightarrow\H^0(\kappa(z),\Z/n)\longrightarrow\H^0(K_{i,z},\Z/n)\oplus\H^0(K_{j,z},\Z/n)\longrightarrow\H^1_z(\O_{C,z}^\h,\Z/n)\to 0
$$
This shows $\H^1_z(\O_{C,z}^\h,\Z/n)\isom\Z/n$,
and since $\H^1_\S(\O_{C,\S},\Z/n)$ is a finite direct sum of these groups,
it is a finite free $\Z/n$-module, of rank $|\S|$.

The result \cite[27.1]{Fch} implies that a free $\Z/n$-submodule
of a $\Z/n$-module is a direct summand.
Therefore we have a decomposition 
$$
\H^0(\kappa(C),\Z/n)\isom\im(\phi_1)\oplus\im(\phi_2)
$$
and since $\H^0(\kappa(C),\Z/n)$ is a finite free $\Z/n$-module,
$\im(\phi_2)$ is a finite free $\Z/n$-module by the structure theorem for finitely generated 
abelian groups.
Similarly, since $\H^1_\S(\O_{C,\S},\Z/n)$ is a finite free $\Z/n$-module,
$$
\H^1_\S(\O_{C,\S},\Z/n)\isom\im(\phi_2)\oplus \cok(\phi_2)
$$
and 
since $\im(\phi_2)$ and $\H^1_\S(\O_{C,\S},\Z/n)$ are finite free $\Z/n$-modules, so is $\cok(\phi_2)$.
Since $\H^1(\O_{C,\S},\Z/n)$ is a $\Z/n$-module, $\cok(\phi_2)$
is a direct summand of $\H^1(\O_{C,\S},\Z/n)$, again by \cite[27.1]{Fch}.
Thus we have a decomposition
$$
\H^1(\O_{C,\S},\Z/n)\isom\cok(\phi_2)\oplus\im(\phi_4)
$$
Now
we set $\Gamma=\im(\phi_4)$, and compute $\rk(\cok(\phi_2))=N+|\S|-m=\beta_C$. 
This proves the $q=1$ part of (b).

The $q=0$ case of (b) is in the proof of Lemma  \ref{gluing}.
Suppose $q=2$.
To show $\H^2(\O_{C,\S},\mu_n)\to\H^2(\kappa(C),\mu_n)$ is injective,
we will show $\H^1(\kappa(C),\mu_n)\to\H^2_\S(\O_{C,\S},\mu_n)$ is onto
and apply the exactness of ($*$).

For each closed point $z\in C_1\cap C_2\subset \S$, we have a diagram
\[\xymatrix{
&\H^1(\kappa(C_1),\mu_n)\oplus\H^1(\kappa(C_2),\mu_n)\ar[r]\ar[d]&\H^2_z(\O_{C,z}^\h,\mu_n)\ar@{=}[d]\\
0\longrightarrow\H^1(\kappa(z),\mu_n)\ar[r]&\H^1(\kappa(C_1)_z,\mu_n)\oplus\H^1(\kappa(C_2)_z,\mu_n)
\ar[r]&\H^2_z(\O_{C,z}^\h,\mu_n)\longrightarrow 0
}\]
We will show that $\H^1(\kappa(C_1),\mu_n)\oplus\H^1(\kappa(C_2),\mu_n)\to\H^2_z(\O_{C,z}^\h,\mu_n)$ is onto,
by showing the downarrow is onto.
Since $z$ is a regular point of each $C_i$, each $\O_{C_i,z}$ is a discrete valuation 
ring with residue field $\kappa(z)$ and fraction field $\kappa(C_i)$,
and we have a diagram of split short exact sequences
\[\xymatrix{
0\ar[r]&\H^1(\O_{C_i,z},\mu_n)\ar[r]\ar[d]&\H^1(\kappa(C_i),\mu_n)
\ar[r]\ar[d]&\H^0(\kappa(z),\Z/n)\ar@{=}[d]\ar[r]&0\\
0\ar[r]&\H^1(\widehat\O_{C_i,z},\mu_n)\ar[r]&\H^1(\kappa(C_i)_z,\mu_n)
\ar[r]&\H^0(\kappa(z),\Z/n)\ar[r]&0\\
}\]
To show the middle downarrow is onto it suffices (by a standard diagram chase) 
to prove that the left downarrow is onto.
Since $\widehat\O_{C_i,z}$ is henselian $\H^1(\widehat\O_{C_i,z},\mu_n)=\H^1(\kappa(z),\mu_n)$,
and by Kummer theory and Hilbert 90 we have $\H^1(\O_{C_i,z},\mu_n)=\O_{C_i,z}^*/n$
and $\H^1(\kappa(z),\mu_n)=\kappa(z)^*/n$.  
Since $\O_{C_i,z}\to \kappa(z)$ is onto and
$\O_{C_i,z}$ is local, the induced map $\O_{C_i,z}^*\to \kappa(z)^*$ is onto, hence $\H^1(\O_{C_i,z},\mu_n)$
maps onto $\H^1(\kappa(z),\mu_n)$.
We conclude $\H^1(\kappa(C_i),\mu_n)\to\H^1(\kappa(C_i)_z,\mu_n)$ is onto. 
Now each map $\H^1(\kappa(C_1),\mu_n)\oplus\H^1(\kappa(C_2),\mu_n)\to\H^2_z(\O_{C,z}^\h,\mu_n)$ is onto. 

Suppose $(b_z)\in\H^2_\S(\O_{C,\S},\mu_n)=\bigoplus_\S\H^2_z(\O_{C,z}^\h,\mu_n)$.
We have just seen that
for each closed point $z\in C_i\cap C_j$ there exists a pair 
$([a_{i,z}t_{i,z}^{e_i}],[a_{j,z}t_{j,z}^{e_j}])\in \kappa(C_i)^*/n\oplus \kappa(C_j)^*/n$
mapping to $b_z$, for $z$-units $a_{k,z}\in\O_{C_k,z}^*$, $z$-uniformizers
$t_{k,z}\in \kappa(C_k)$, and integers $e_k$, for $k=i,j$.
Let $v_{k,z}$ be the discrete valuation on $\kappa(C_k)$ determined by $z$.
By standard approximation 
(e.g. \cite[XII.1.2]{Lang}) there exist elements $a_k,t_k\in \kappa(C_k)$ such that 
$$
v_{k,z}(a_k-a_{k,z})>0\quad\text{and}\quad v_{k,z}(t_k-t_{k,z})>1\quad\text{for all $z$.}
$$
The image of $a_kt_k^{e_k}$ in $\kappa(C_k)_z^*/n$ is $[a_{k,z}t_{k,z}^{e_k}]$.
Therefore the $m$-tuple 
$$
([a_k t_k^{e_k}])\in\H^1(\kappa(C),\mu_n)
$$
maps to $(b_z)$.
This proves 
%
%
the induced map $\H^1(\kappa(C),\mu_n)\to\H^2_\S(\O_{C,\S},\mu_n)$ is onto, and completes the proof.

\hfill $\blacksquare$

We will soon need the following technical lemma in order to replace $X_0$ with $C$.

\Lemma\label{reduce}
Suppose $A$ is a noetherian ring.
Then $(\Frac A)_\red = \Frac(A_\red)$ if and only if A has no embedded primes.
\rm

\begin{proof}
By definition $\Frac A=S^{-1}A$, where $S=A-\bigcup_{\Ass A}\frak p$, 
and $S^{-1}N_A=N_{S^{-1}A}$ by \cite[3.12]{AM}, hence $S^{-1}(A_\red)=(S^{-1}A)_\red$.  
It remains to show that
$S^{-1}(A_\red) = \Frac(A_\red)$ if and only if $A$ has no embedded primes.  
By $S^{-1}(A_\red)$ of course we mean $f(S)^{-1}(A_\red)$, where $f:A\to A_\red$.  
This localization equals the localization with respect to the multiplicative set T, 
where T is the saturation of $f(S)$ in $A_\red$, and this is the complement of the union of 
prime ideals of $A_\red$ that don't meet $f(S)$ by \cite[Exercise 3.7]{AM}.  
Thus $S^{-1}(A_\red)=\Frac(A_\red)$ if and only if the union of the primes of $A_\red$ that don't meet 
$f(S)$ equals the union of the associated primes of $A_\red$, which are just the minimal primes of 
$A_\red$.  But $A$ and $A_\red$ have identical underlying topological spaces, and the primes of 
$A_\red$ that don't meet $f(S)$ correspond to the primes of $A$ that don't meet $S$, i.e., 
the associated primes. These correspond to the minimal primes of $A_\red$ if and only if 
the associated primes of $A$ are the minimal primes of $A$, i.e., 
$A$ has no embedded primes.
\end{proof}

\Theorem\label{map}
Assume the setup of (\ref{setup}), with $X$ connected.
Then for $q\geq 0$
there is a map $$\lambda:\H^q(\O_{C,\S},\Lambda)\to\H^q(K(X),\Lambda)$$ 
and a commutative diagram
\begin{equation}\label{lambdadiagram}
\xymatrix{
\H^q(\O_{C,\S},\Lambda)\ar[r]^-\lambda\ar[d]_{\oplus\res_i}&\H^q(K(X),\Lambda)\ar[d]^{\oplus_i\res_i}\\
\H^q(\kappa(C),\Lambda)\ar[r]^-\inf&\bigoplus_i\H^q(K(X)_{C_i},\Lambda)
}\end{equation}
such that if $\alpha_0\in\H^q(\O_{C,\S},\Lambda)$ and $\alpha=\lambda(\alpha_0)$
then:
\begin{enumerate}
\item[a)]
$\alpha$ is defined at the generic points of $C_i$, and $\alpha(C_i)=\res_i(\alpha_0)$.
\item[b)]
The ramification locus of $\alpha$ (on $X$) is contained in $\mathscr D_\S$.
\item[c)]
If $D\in\mathscr D_\S$ is prime and $z=D\cap C$,
then $\partial_D\cdot\lambda=\inf_{\kappa(z)|\kappa(D)}\cdot\partial_z$.
\item[d)]
If $\alpha_0$ is unramified at a closed point $z$,
and $D$ is any (horizontal) prime lying over $z$,
then $\alpha$ is unramified at $D$, and has value
$\alpha(D)=\inf_{\kappa(z)|\kappa(D)}(\alpha_0(z))$.
\end{enumerate}
\rm

\Pf
Let $D_0$ be an effective divisor on $C$ that avoids $\S$, 
let $D\in\mathscr D_\S$ be the distinguished lift of $D_0$, set $U=X-D$,
and set $U_0=C-D_0$.
Since $X$ and $D$ are regular and $D$ has pure codimension $1$, we have 
$\H^0(X)\isom\H^0(U)$, and an exact Gysin sequence
$$
0\longrightarrow\H^1(X)\longrightarrow\H^1(U)
\xrightarrow{\;\partial_D\;}\H^0(D,-1)\longrightarrow\H^2(X)\longrightarrow\cdots
$$
by Gabber's absolute purity theorem (\cite[Theorem 2.1.1]{Fuj02}) and the standard
construction of the Gysin sequence (\cite[Section 3.2]{CT95}).
(Note that the result in \cite{Fuj02} is stated for the $\Lambda=\Z/n$ case only,
but the result holds in general since the sheaves 
$\mathscr H_D^q(X)$ and $\mathscr H_D^q(X,\Z/n)$ are locally isomorphic,
and the morphism $i^*\Lambda(-1)\to\mathscr H_D^2(X)$ is canonical.)
We use the notation $\partial_D$ since this map
is compatible with the one defined above on $\H^q(K(X))$ when $D$ is prime.

We may replace $X_0$ by $C=X_{0,\red}$ in the cohomological computations below
since $\Lambda$ is finite and $n$ is prime-to-$p$,
by \cite[V.2.4(c)]{M} (see also \cite[II.3.11]{M}).  To substitute $\O_{C,\S}$
and $\kappa(C)$ for $\O_{X_0,\S}$ and $\kappa(X_0)$ we must check that the former
are the canonical reduced quotients of the latter.
But the ring $\O_{X_0,\S}$ can by obtained by localizing some affine open subset $\Spec A_0$
containing $\S$ (which exists since $X_0/k$ is projective)
with respect to the multiplicative set $T=A_0-\bigcup_\S\frak m_z$.
Since $\O_{C,\S}$ is obtained by localizing $A_{0,\red}$ with respect to the image of $T$ in $A_{0,\red}$,
we have $\O_{C,\S}=(\O_{X_0,\S})_\red$ since the formation of the nilradical commutes
with localization (see e.g. \cite[3.12]{AM}).

To show $\kappa(C)=\kappa(X_0)_\red$ it suffices to show $X_0$ has no embedded points
by Lemma~\ref{reduce}.
But if $z$ is any closed point of $X$ then $\O_{X,z}$ is a regular local ring,
and a local equation for the closed fiber $\O_{X,z}\otimes_R k$ passing through $z$ is given by
the uniformizer $p$ in $R$.  Since $\O_{X,z}$
is factorial and at most two components of $X_0$ pass through $z$ we have
$p=\pi_1^{e_1}\pi_2^{e_2}$ for primes $\pi_i$ and numbers $e_i\geq 0$.
The associated primes of $\O_{X,z}/(\pi_1^{e_1}\pi_2^{e_2})$ are evidently just
the $(\pi_i)$, which shows $X_0$ has no embedded point at $z$.

Since $D_0$ is a disjoint union of regular closed points,
by ($**$) and the work that immediately precedes it we have
$\H^0(C)\isom\H^0(U_0)$ and an exact sequence
$$
0\longrightarrow\H^1(C)\longrightarrow\H^1(U_0)\xrightarrow{\;\partial_{D_0}\;}\H^0(\kappa(D_0),-1)
\longrightarrow \H^2(C)\longrightarrow\cdots
$$
Thus we have a commutative ladder
\[\xymatrix{
0\ar[r]&\H^1(X)\ar[r]\ar[d]&\H^1(U)\ar[r]^-{\partial_D}\ar[d]
&\H^0(D,-1)\ar[r]\ar[d]&\H^2(X)\ar[r]\ar[d]&\cdots\\
0\ar[r]&\H^1(C)\ar[r]&\H^1(U_0)\ar[r]^-{\partial_{D_0}}&\H^0(D_0,-1)
\ar[r]&\H^2(C)\ar[r]&\cdots
}\]
Since $R$ is complete,
$\H^q(X)\to\H^q(C)$ and $\H^q(D,-1)\to\H^q(D_0,-1)$ are
isomorphisms for $q\geq 0$ by proper base change (\cite[VI.Corollary 2.7]{M}).
Therefore, in light of the isomorphisms in degree zero and the 5-lemma in degree $q\geq 1$, 
we obtain isomorphisms
$$\H^q(U)\isim\H^q(U_0)$$
for $q\geq 0$.
Let $\widetilde U$ denote the inverse limit over all these open sets $U$ (this is a scheme by \cite[8.2.3]{EGAIV:c}).
Then $\H^q(\widetilde U)$ is the direct limit of the $\H^q(U)$ by \cite[III.1.16]{M}, and since
the direct limit functor is exact we have an isomorphism
$\H^q(\widetilde U)\isim\H^q(\O_{C,\S})$.
Composing the inverse with $\H^q(\widetilde U)\to\H^q(K(X))$ yields our lift
$$\lambda:\H^q(\O_{C,\S})\longrightarrow\H^q(K(X))$$
The commutative diagram (\ref{lambdadiagram}) follows by applying cohomology to the diagram
\[\xymatrix{
\Spec\O_{C,\S}\ar[r]&\widetilde U&\Spec K(X)\ar[l]\\
\Spec \kappa(C_i)\ar[r]\ar[u]&\Spec\O_{K(X)_{C_i}}\ar[u]&\Spec K(X)_{C_i}\ar[u]\ar[l]
}\]
incorporating the isomorphisms induced by the upper and lower left horizontal arrows.
If $\alpha_0\in\H^q(\O_{C,\S})$ and $\alpha=\lambda(\alpha)$,
then since $\widetilde U$ contains $\S$ and the generic points
of the $C_i$, $\alpha$ is defined at these points,
and the formula $\alpha(C_i)=\res_{\O_{C,\S}|\kappa(C_i)}(\alpha_0)$
follows immediately from (\ref{lambdadiagram}), proving (a).
 
If $D$ is a horizontal prime divisor not in $\mathscr D_\S$, then the generic point
$\Spec\kappa(D)$ is contained in $\widetilde U$, hence the map $\H^q(\widetilde U)\to\H^q(K(X)_D)$
factors through $\H^q(\O_{K(X)_D})$, which shows $\partial_D\cdot\lambda=0$.
Thus the ramification locus of any element in the image of $\lambda$ must be contained
in $\mathscr D_\S$, proving (b).
Now if $D\in\mathscr D_\S$ is prime and $z=D\cap C$ then 
$D$ is the prime spectrum of a complete local ring with residue field $\kappa(z)$, and the isomorphism
$$\H^{q-1}(D,-1)\isim\H^{q-1}(z,-1)=\H^{q-1}(\kappa(z),-1)$$
is the standard identification.
Thus the formula $\partial_D\cdot\lambda=\inf_{\kappa(z)|\kappa(D)}\cdot\partial_{z}$ is immediate
by the commutative ladder of Gysin sequences above, proving (c).

Suppose $\alpha=\lambda(\alpha_0)$ has ramification locus $D_\alpha$,
then $D_\alpha\in\mathscr D_\S$.  Set $U=X-D_\alpha$.
If $\alpha_0$ is unramified at a point $z$, then $\alpha$ is unramified at 
every prime divisor $D$ lying over $z$.
For if $D\in\mathscr D_\S$ then $\partial_D(\alpha)=\inf(\partial_{z}(\alpha_0))$
by the formula just proved,
and if $D\not\in\mathscr D_\S$ then $\partial_D(\alpha)=0$ since $D_\alpha\in\mathscr D_\S$.
Thus if $\alpha_0$ is unramified at $z$, and $D$ is a prime divisor lying over $z$,
then $U$ contains $D$.
%
The maps $z=\Spec \kappa(z)\to U_0$ and $D\to U$ then induce a commutative diagram
\[\xymatrix{
\H^q(U)\ar[r]^-{\res}\ar[d]^{\res}&\H^q(D)\ar[r]^-{\res}\ar[d]^{\res}&\H^q(\kappa(D))\\
{\H^q(U_0)}{\ar[r]^-{\res}}&{\H^q(z)}\ar[ur]_{\inf}}
\]
Both vertical down-arrows are isomorphisms by proper base change. 
The inverse of the left one is $\lambda$ by definition,
and the composition of the inverse of the right one and the restriction
$\H^q(D)\to\H^q(\kappa(D))$ is inflation, as shown. 
Since $\kappa(D)$ is complete, the top composition of horizontal restrictions factors through the restriction
$\H^q(U)\to\H^q(\widehat\O_{X,D})$ and 
the bottom factors through the restriction
$\H^q(U_0)\to\H^q(\widehat\O_{C,z})$. 
Since these are restriction maps, the images of $\alpha$ and $\alpha_0$ are 
the values $\alpha(D)$ and $\alpha_0(z)$.
We conclude $\inf_{\kappa(z)|\kappa(D)}(\alpha_0(z))=\alpha(D)$, as in (d).

\hfill $\blacksquare$

\Paragraph\label{pi}
By weak approximation \cite[Lemma]{Sa98}
there exists a $\pi\in K(X)$ such that 
$$
\div\,\pi=C+E
$$ 
where $E$ contains no components of $C$, and avoids any finite set of points $\mathcal N$.
We fix such a $\pi$ for $\mathcal N$ containing $\S$. 
For each $i$, the choice of $\pi$ determines a noncanonical ``Witt'' isomorphism 
\[
\H^q(\kappa(C_i))\oplus\H^{q-1}(\kappa(C_i),-1)\isim\;\H^q(K(X)_{C_i})
\]
Taking $(\alpha,\theta)$
to $\alpha+(\pi)\cdot\theta$, where $\alpha$ and $\theta$ are inflated from 
$\kappa(C_i)$ to $K(X)_{C_i}$,
$(\pi)$ is the image of $\pi$ in $\H^1(K(X)_{C_i},\mu_n)$, and $(\pi)\cdot\theta$
is the cup product.
Although we cannot in general lift all of
$\bigoplus_i\H^q(K(X)_{C_i})$ to $\H^q(K(X)$, we can now prove the following.

\Corollary\label{splits}
Let $(\pi)$ denote the image of $\pi$ in $\H^1(K(X),\mu_n)$.
The choice of $\mathscr D_\S$ and $\pi$ determines a homomorphism for $q\geq 1$,
\begin{align*}
\lambda:\H^q(\O_{C,\S},\Lambda)\oplus\H^{q-1}(\O_{C,\S},\Lambda(-1))&\longrightarrow\;\H^q(K(X),\Lambda)\\
(\alpha_0,\theta_0)&\longmapsto \lambda(\alpha_0)+(\pi)\cdot\lambda(\theta_0)
\end{align*}
such that $(\bigoplus_i\res_{K(X)|K(X)_{C_i}})\cdot\lambda
=\bigoplus_i(\inf_{\kappa(C_i)|K(X)_{C_i}}\cdot\res_{\O_{C,\S}|\kappa(C_i)})$.
\rm

\Pf
This is an immediate consequence of Theorem \ref{map}.

\hfill $\blacksquare$

\Remark\label{connected}
a) Theorem \ref{map} and Corollary \ref{splits} 
apply with obvious amendments to the case where $X$ is not connected.  
For if $X=\coprod_k X_k$ is a decomposition into connected components,
then $K(X)=\prod_k K(X_k)$, $X_0=\coprod_k(X_k)_0$, $\O_{X_0,\S}=\prod_k\O_{(X_k)_0,\S_k}$
(where $\S_k=\S\cap X_k$),
all of the cohomology groups break up into direct sums, and we define the map $\lambda$ 
to be the direct sum of the maps on the summands.
This will come up in the next section.

b) If $X$ is smooth, then $\S$ is empty, and $\O_{C,\S}=\kappa(C)$.
By Witt's theorem we have 
$\H^q(K(X)_C,\Lambda)\isom\H^q(\kappa(C),\Lambda)\oplus\H^{q-1}(\kappa(C),\Lambda(-1))$,
and we obtain a map
$$\lambda:\H^q(K(X)_C,\Lambda)\longrightarrow\H^q(K(X),\Lambda)$$
that splits the restriction map.
This is the map of \cite{BMT}.

\Paragraph{Completely split characters.}
In \cite[2.1]{Sai85} Saito defines a completely split covering
of a noetherian scheme $X$ to be a finite \'etale cover $Y\to X$ such that 
$Y\times_X \Spec \kappa(x)=\coprod\Spec\kappa(x)$, for all closed points $x\in X$. 
We abuse Saito's terminology (see Remark\eqref{remark} below) and in the setup of \eqref{setup} denote by
$\H^1_\cs(C,\Z/n)$ the kernel of the map $\H^1(\O_{C,\S},\Z/n)\to\H^1(\kappa(C),\Z/n)$
in Lemma~\ref{injects}.
If $\beta\in\H^1_\cs(C,\Z/n)$ then $\partial_z(\beta)=0$ for all closed points $z\in C-\S$ since
$\partial_z$ factors through $\kappa(C_i)_z$.
Therefore $\beta$ is defined on $C$, hence $\H^1_\cs(C,\Z/n)\leq\H^1(C,\Z/n)$.
Let $\H^1_\cs(X,\Z/n)$ denote the preimage of $\H^1_\cs(C,\Z/n)$ under the proper
base change isomorphism.

\Proposition\label{cs}
Assume the setup of (\ref{setup}).
Then elements of $\H^1_\cs(C,\Z/n)$ are trivial at all points of $C$,
and the nontrivial elements of $\H^1_\cs(X,\Z/n)$ are trivial at all points of $X$ except for the
generic point $\Spec K(X)$, where they are nontrivial.
\rm

\begin{proof}
Suppose $\beta_0\in\H^1_\cs(C)$.
Then $\beta_0$ is trivial at each generic point of $C$ by definition of $\H^1_\cs(C)$.
If $z\in C$ is a closed point lying on the irreducible component $C_i$ then
the map $\H^1(C)\to\H^1(\kappa(z))$ factors through $\H^1(C_i)$.
Since $C_i$ is regular the map $\H^1(C_i)\to\H^1(\kappa(C_i))$ is injective by purity, 
and consequently $\beta_0(z)=0$ by definition.
Thus the elements of $\H^1_\cs(C)$ are trivial at all points of $C$.

Suppose $\beta=\lambda(\beta_0)\in\H^1_\cs(X)$.
If $x\in X$ is a generic point of some irreducible component
$C_i$ of $C$ then the image of $\beta$ in $\H^1(\kappa(C_i))$
is zero since the map $\H^1_\cs(X,\Z/n)\to\H^1(\kappa(C_i))$ factors through $\H^1_\cs(C)$.
If $x$ is the generic point of a horizontal divisor $D$ with closed point $z$ then 
$\beta(D)=\inf_{\kappa(z)|\kappa(D)}(\beta_0(z))$ by Theorem~\ref{map}(d),
and this is zero since $\beta_0(z)=0$.
If $z$ is a closed point of $X$ then $z$ is on $C$,
and the map $\H^1_\cs(X)\to\H^1(\kappa(z))$ factors through $\H^1_\cs(C)$, hence $\beta$
is trivial at $z$.
Finally, since $X$ is regular the map $\H^1(X)\to\H^1(K(X))$ is injective by purity,
hence $\beta$ is nontrivial at the generic point of $X$.
\end{proof}

\Remark\label{remark}
Proposition~\ref{cs} shows the elements of $\H^1_\cs(C,\Z/n)$ are completely split in the sense of \cite{Sai85}.
However, in the general case our $\H^1_\cs(C,\Z/n)$ does not account for elements 
that are split at every closed point but nontrivial at generic points
of $C$.  This is not an issue if $k$ is finite as shown by Saito in \cite[Theorem 2.4]{Sai85},
since then the $C_i$ have no nontrivial completely split covers, essentially
by Cebotarev's density theorem (see \cite[Lemma 1.7]{Ras}).


\section{Index Calculation in the Brauer Group}

\Paragraph{Cyclic Covers.}
If $U$ is any scheme, and $\bar u$ is a geometric point,
the fiber functor defines a category equivalence between (finite) \'etale covers of $U$ 
and finite continuous $\pi_1(U,\bar u)$-sets, yielding a canonical isomorphism
$$
\H^1(U,\Z/n)\isom\H^1(\pi_1(U,\bar u),\Z/n)=\Hom_\cont(\pi_1(U,\bar u),\Z/n)
$$
(see \cite[I.2.11]{FK}).
If $\theta\in\H^1(U,\Z/n)$, we will write $U[\theta]$ for the finite cyclic \'etale 
cover determined by $\theta$.  If $U=\Spec A$ is affine, we will write $A[\theta]$
for the corresponding ring, or $A(\theta)$ if $A$ is a field.
If $U$ is a connected normal scheme, and $\theta\in\H^1(U,\Z/n)$ has order $m$, 
then $U[\theta]$ is a disjoint sum of $n/m$ connected $\Z/m$-Galois covers of $U$.

\Lemma\label{cycliccovers}
Assume the setup of (\ref{setup}).
Let $\theta_0\in\H^1(\O_{C,\S},\Z/n)$ be a (tamely ramified) character
with ramification divisor $D_0$ on $C$.
Then the (tame) ramification divisor of $\theta=\lambda(\theta_0)$ is
the distinguished lift $D\in\mathscr D_\S$ of $D_0$ on $X$, and 
$\theta$ defines a tamely ramified cover $\rho:Y\to(X,D)$
as in Lemma~\ref{covers}.
Restriction to $C$ yields a tamely ramified cover
$\rho_0:C_Y\to(C,D_0)$ such that $\O_{C_Y,\S_Y}=\O_{C,\S}[\theta_0]$, 
and the reduced closed fiber
$C_Y$ of $Y$ is the normalization of $C$ in $\kappa(C)(\theta_0)$.
\rm

\Pf
The lift $\theta$ is tamely ramified with respect to $D$ by Theorem \ref{map}.
Let $Y$ be the normalization of $X$ in $L=K(X)(\theta)$.
Since $D$ is in $\mathscr D_\S$ and $X/R$ satisfies the setup of \eqref{setup}, by 
Lemma~\ref{covers} $\rho:Y\to(X,D)$ is a tamely ramified cover,
$Y/R$ is a regular relative curve with reduced closed fiber $C_Y$,
the irreducible components of $C_Y$ are regular with singular points $\S_Y$,
and $D_Y$ is in $\mathscr D_{Y,\S_Y}$.

Let $U=X-D$, $V=U\times_X Y$, $U_0=U\times_X X_0$ and $V_0=V\times_Y Y_0$.
Then $\theta_0\in\H^1(U_0)$.
We have $\H^1(U)\isom\Hom(\pi_1(U),\Z/n)$ and $\H^1(U_0)\isom\Hom(\pi_1(U_0),\Z/n)$,
and the restriction map $\res:\H^1(U)\to\H^1(U_0)$, which sends $\theta$ to $\theta_0$, 
is induced by the natural map $\pi_1(U_0)\to\pi_1(U)$,
which is induced on covers by $W\mapsto W\times_U U_0$.
Therefore $V_0=U_0[\theta_0]$.

We show $\O_{C_Y,\mathcal S_Y}=\O_{C,\S}[\theta_0]$.
If $U_0=\Spec A_0\subset C-D_0$ is a dense affine open subset of $C$ containing $\S$,
then its preimage in $C_Y$ is a dense affine open subset $V_0=\Spec B_0$ containing $\mathcal S_Y$.
By base change we have $B_0=A_0[\theta_0]$, and $S^{-1}B_0=\O_{C,\S}[\theta_0]$,
where $S=A_0-\left(\bigcup_{x\in\S}\frak m_x\right)$ is the multiplicative set defining $\O_{C,\S}$.
Since $\mathcal S_Y=\rho^{-1}\S$,
the saturation $T$ of $S$ in $B_0$ is $T=B_0-\left(\bigcup_{y\in\mathcal S_Y}\frak m_y\right)$,
which shows $S^{-1}B_0=\O_{C_Y,\mathcal S_Y}$.
Therefore $\O_{C_Y,\mathcal S_Y}=\O_{C,\S}[\theta_0]$, and it follows immediately
that $\O_{C_Y,\S_Y}=\O_{C,\S}[\theta_0]$.

The map $C_Y\to C$ is finite by base change,
and since each irreducible component of $C_Y$ is regular by Lemma~\ref{structure} 
each irreducible component of $C_Y$ is the normalization of a component of $C$ in 
a field extension which is a direct factor of $\kappa(C_Y)$.
Equivalently, $C_Y$ is the normalization of $C$ in $\kappa(C)(\theta_0)$, by \cite[6.3.7]{EGAII}.
This completes the proof.

\hfill $\blacksquare$

\Paragraph{Index.}\label{index}
Assume the setup of \eqref{setup} with $R=\Z_p$ and $\Lambda=\mu_n$.
Fix $\alpha_C\in\H^2(\O_{C,\S})$ and $\theta_C\in \Gamma\leq\H^1(\O_{C,\S},-1)$, 
and write
\begin{align*}
\alpha_C&=(\alpha_{C_1},\dots,\alpha_{C_m})\in\H^2(\kappa(C))\\
\theta_C&=(\theta_{C_1},\dots,\theta_{C_m})\in\H^1(\kappa(C),-1)
\end{align*}
as per Lemma~\ref{injects}(b).
Suppressing the inflation maps, we form the element
$$
\gamma_C=
\alpha_C+(\pi)\cdot\theta_C=(\alpha_{C_1}+(\pi)\cdot\theta_{C_1},
\dots,\alpha_{C_m}+(\pi)\cdot\theta_{C_m})
\in\bigoplus_{i=1}^m\H^2(K(X)_{C_i})
$$
where $\pi$ is as in \eqref{pi}.
Let $\theta=\lambda(\theta_C)$, and let $Y$ be the normalization of $X$ in $K(X)(\theta)$,
as in Lemma~\ref{cycliccovers}.
Then $C_Y$ is the reduced closed fiber of $Y$, and we write $C_{i,Y}$ for the preimage of $C_i$,
so that $\kappa(C_{i,Y})=\kappa(C_i)(\theta_{C_i})$.
Thus
$$
\alpha_{C_Y}=(\alpha_{C_{1,Y}},\dots,\alpha_{C_{m,Y}})
\;\in\;\bigoplus_{i=1}^m\H^2(\kappa(C_{i,Y}))
$$
Note $\kappa(C_{i,Y})$ is a product of the function fields of the irreducible components of $C_{i,Y}$.
By the (well-known) Nakayama-Witt index formula,
$$
\ind(\alpha_{C_i}+(\pi)\cdot\theta_{C_i})=|\theta_{C_i}|\cdot\ind(\alpha_{C_{i,Y}})
$$
We have $|\theta_C|=\lcm_i\{|\theta_{C_i}|\}$ by Lemma~\ref{injects}(b). 
We now {\it define}
\begin{align*}
\ind(\alpha_{C_Y})\;&\df\;\lcm_i\{\ind(\alpha_{C_{i,Y}})\}\\
\ind(\gamma_C)\;&\df\;|\theta_C|\cdot\ind(\alpha_{C_Y})
\end{align*}

\Theorem\label{preservesindex}
Assume the setup of (\ref{setup}) with $R=\Z_p$.
Let $\Gamma\leq\H^1(\O_{C,\S})$ be as in Lemma~\ref{injects}(b).
Then the map $\lambda:\H^2(\O_{C,\S})\oplus \Gamma\to\H^2(K(X))$
preserves index.
\rm

\Pf
We may assume $X$ is connected.
We identify $\Gamma$ with the image of $\H^1(\O_{C,\S},-1)$ in $\H^1(\kappa(C),-1)$, 
as in Lemma~\ref{injects}(b),
and adopt the notation of (\ref{index}).
Set $\gamma=\lambda(\gamma_C)$, $\alpha=\gamma(\alpha_C)$, and $\theta=\lambda(\theta_C)$,
so that $\gamma=\alpha+(\pi)\cdot\theta$ as in Corollary \ref{splits}.
By restricting to connected components if necessary
we may assume that $Y$ is connected, hence that $K(Y)$ is a field.
Even so, the cyclic-Galois \'etale $\kappa(C_i)$-algebra $\kappa(C_{i,Y})$ may not be a field.
Since $\theta_C$ is in $\H^1(\O_{C,\S},-1)$ the ramification divisor $D_0$ of $\theta_C$ avoids $\S$,
and the distinguished lift
$D\in\mathscr D_\S$ of $D_0$ is the ramification divisor of $\theta$ by Theorem \ref{map}.
By Lemma~\ref{cycliccovers} $\rho:Y\to (X,D)$ is a cyclic tamely ramified cover,
and by Lemma~\ref{structure} $Y/R$ satisfies the properies of \eqref{setup},
with reduced closed fiber $C_Y$, $\S_Y=\rho^{-1}\S$ the singular points of $C_Y$,
and $\mathscr D_Y$ generated by $D_{Y,\red}$ and the preimages of the other distinguished
divisors of $X$.

The index of $\gamma$ cannot exceed $|\theta|\ind(\alpha_Y)$.
For if $M/K(Y)$ is a separable maximal subfield of the division algebra
associated with $\alpha_Y$, then $M/K(X)$ splits $\gamma$,
and has degree $|\theta|\ind(\alpha_Y)$.
Since in our case $|\theta|=[K(Y):K(X)]=[\kappa(C_Y):\kappa(C)]=|\theta_C|$, 
to prove the theorem it is enough
to prove $\ind(\alpha_Y)=\ind(\alpha_{C_Y})$.


By Lemma~\ref{cycliccovers} $\O_{C,\S}[\theta_C]=\O_{C_Y,\mathcal S_Y}$.
Each $\kappa(C_{i,Y})=\kappa(C_i)(\theta_{C_i})$ is a product of global fields, 
and by class field theory
the division algebra associated with the restriction of
$\alpha_{C_{i,Y}}$ to each field component is cyclic.
Since $\alpha_{C_Y}$ is in $\H^2(\O_{C_Y,\S_Y})$ (by Lemma~\ref{cycliccovers} and Lemma~\ref{gluing}),
$\alpha_{C_{i,Y}}$ is unramified at $\mathcal S_Y$.
By Grunwald-Wang's theorem there exists a tuple
$\psi_{C_Y}=(\psi_{C_{1,Y}},\dots,\psi_{C_m,Y})\in\bigoplus_i\H^1(\kappa(C_{i,Y}),-1)$
such that $|\psi_{C_{i,Y}}|=\ind(\alpha_{C_{i,Y}})$,
$\kappa(\psi_{C_{i,Y}})$ splits $\alpha_{C_{i,Y}}$, and such that
the $\psi_{C_{i,Y}}$ are unramified and equal at the local fields defined by the singular points $\mathcal S_Y$.
Then $\psi_{C_Y}$ comes from $\H^1(\O_{C_Y,\mathcal S_Y},-1)$ by Lemma~\ref{gluing},
and $|\psi_{C_Y}|=\ind(\alpha_{C_Y})$.

By Theorem \ref{map} (and Remark \ref{connected}(a) if $Y$ is not connected)
we have a map $\lambda_Y:\H^1(\O_{C_Y,\mathcal S_Y},-1)\to\H^1(K(Y),-1)$.
Since the distinguished divisors on $Y$ are the (reduced) preimages of those on $X$,
$\lambda_Y$ is compatible with $\lambda$ and the residue maps.
Set $\psi=\lambda_Y(\psi_{C_Y})$, and let $D_\psi$ 
denote the distinguished lift of the ramification divisor $D_{\psi_{C_Y}}$ of $\psi_{C_Y}$ on $C_Y$.
By Lemma~\ref{cycliccovers}, $\psi$ determines a cyclic tamely ramified cover $\sigma:Z\to(Y,D_\psi)$
with reduced closed fiber $C_Z$, 
such that $C_Z$ is cyclic and tamely ramified over $(C_Y,D_{\psi_{C_Y}})$,
inducing $\psi_{C_Y}$.
Since $\kappa(C_{i,Z})$ splits $\alpha_{C_{i,Y}}$, $\alpha_{C_Z}=0$.

Again we may assume $Z$ is connected.
By construction, $[K(Z):K(Y)]=|\psi|=|\psi_{C_Y}|=[\kappa(C_Z):\kappa(C_Y)]=\ind(\alpha_{C_Y})$.
By (\ref{lambdadiagram}) we have 
$\ind(\alpha_{C_Y})\leq \ind(\alpha_Y)$, and it remains to show $\alpha_Z=0$.
It is then enough to show $\alpha_Z$ is unramified with respect to all Weil divisors on $Z$,
by \cite[Lemma 3.5]{BMT}.

Let $D'=D\cup\rho(D_\psi)$.  
Since $D_\psi\in\mathscr D_{Y,\S_Y}$, $\rho(D_\psi)\in\mathscr D_\S$, and the composition
$\rho':Z\to (X,D')$ is a tamely ramified cover.
Since $X$ is regular, $\rho'$ is (finite and) flat by Lemma~\ref{structure}, 
and so the image of any prime divisor $J$ of $Z$
is a prime divisor $\rho'(J)=I$ of $X$.
By the functoriality of the residue maps, $\alpha_Z$ can only be ramified at prime divisors lying
over irreducible components of $D_\alpha$.
By Theorem~\ref{map}(a) $D_\alpha$ is in $\mathscr D_\S$, and the divisors of $Z$ lying over $D_\alpha$
are distinguished by Lemma~\ref{structure}(c).
Thus it is enough to show that $\alpha_Z$ is unramified at these distinguished divisors.
Clearly we may assume $D_\alpha$ is irreducible.

Let $E\subset Z$ be a (distinguished) prime divisor lying over $D_\alpha\in\mathscr D_\S$.
Since $D_\alpha\cap D'=\varnothing$ or $D_\alpha\subset D'$,
we have $e(v_E/v_{D_\alpha})=e(v_{E_0}/v_{D_{\alpha_C}})=e$ for some $e\geq 1$,
by Lemma~\ref{lemma2}.
By Lemma~\ref{map} and the functorial behavior of the residue and restriction
maps, we have a commutative diagram
\[\xymatrix{
\H^2(\O_{C,\S})\ar[d]_{\partial_{D_{\alpha_C}}}\ar[r]^-\lambda\ar@{}[dr]|{1}&\H^2(K(X))\ar[d]_{\partial_{D_\alpha}}
\ar[r]^-\res\ar@{}[dr]|{2}&\H^2(K(Z))\ar[d]^{\partial_E}\\
\H^1(\kappa(D_{\alpha_C}),-1)\ar[r]&\H^1(\kappa({D_\alpha}),-1)\ar[r]^-{e\cdot\res}\ar@{}[dr]|{3}&\H^1(\kappa(E),-1)\\
&\H^1(\kappa(D_{\alpha_C}),-1)\ar[r]^-{e\cdot\res}\ar[u]^\inf\ar@{}[dr]|{4}&\H^1(\kappa(E_0),-1)\ar[u]_\inf\\
&\H^2(\O_{C,\S})\ar[u]^{\partial_{D_{\alpha_C}}}\ar[r]^\res&\H^2(\O_{C_Z,\mathcal S_Z})\ar[u]_{\partial_{E_0}}\\
}\]
Since $\alpha=\lambda(\alpha_C)$, 
$\partial_E(\alpha_Z)=e\cdot(\partial_{D_{\alpha_C}}(\alpha_C))_{\kappa(E)}$
by squares (1) and (2), and by square (4),
$\partial_{E_0}(\alpha_{C_Z})=e\cdot(\partial_{D_{\alpha_C}}(\alpha_C))_{\kappa(E_0)}$.
Therefore $\partial_E(\alpha_Z)=\inf_{\kappa(E_0)|\kappa(E)}(\partial_{E_0}(\alpha_{C_Z})$ 
by square (3).
Since 
$\alpha_{C_Z}=0$, we conclude $\partial_E(\alpha_Z)=0$,
as desired.
This proves the theorem.

\hfill $\blacksquare$

\section{Noncrossed Products and Indecomposable Division Algebras}

A (finite-dimensional) 
division algebra $D$ central over a field $F$ is called a {\it noncrossed product}
if it has no Galois maximal subfield.  Its algebra structure then cannot be given
by a Galois 2-cocycle, counter to almost all known division algebra constructions
(see \cite{Ha85} for a construction of a noncrossed product).  
Noncrossed product division algebras were long conjectured to be fictional, until they
were shown to exist by Amitsur in \cite{Am72}. 

We say $D$ is {\it indecomposable} if it does not contain a subalgebra that is also
central over $F$, or equivalently if it is not an $F$-tensor product of two nontrivial
$F$-division algebras.  It is not hard to show that all division algebras of composite period
are decomposable, and that all division algebras of equal prime-power period 
and index are indecomposable, but it is nontrivial to construct indecomposable division algebras of unequal
prime-power period and index.  
The first examples appeared in \cite{ART} and in \cite{Sa79}.
For additional discussion of either of these
topics, see almost any survey treating open problems on division algebras, such as
\cite{ABGV}, \cite{Am82}, or \cite{Sa92}.

We can use Theorem \ref{preservesindex} to prove the existence of 
noncrossed product and indecomposable division algebras over the function
field $F$ of any $p$-adic curve $X_{\Q_p}$.
Noncrossed products over $K(t)$ for $K$ a local field were first
constructed in \cite{Br01a}, and then constructed more sytematically over
the function field of a smooth relative $\Z_p$-curve in \cite{BMT}.
Indecomposable division algebras of unequal period and index were also constructed in 
\cite{BMT}, over the same types of fields.
Modulo gluing,
the method we use below is the same as the one
used in \cite[Theorem 4.3, Corollary 4.8]{BMT}.

\Theorem\label{noncrossed}
Let $F/\Q_p$ be a finitely generated field extension of transcendence degree one.
Let $X/\Z_p$ be a regular relative curve with function field
$F$, let $C_i$ be a reduced irreducible
component of the closed fiber, let $\ell\neq p$ be a prime, 
and let $r$ and $s$ be numbers that are maximal such that $\mu_{\ell^r}\subset \kappa(C_1)$
and $\mu_{\ell^s}\subset \kappa(C_1)(\mu_{\ell^{r+1}})$.
Then there exist noncrossed product $F$-division algebras of period and index
as low as
$\ell^{s+1}$ if $r=0$, and $\ell^{2r+1}$ if $r\neq 0$.
\rm

\Pf
We may assume (without changing $r$ and $s$) that $C$ has regular irreducible components,
at most two of which meet at any closed point of $X$.
The idea is to use the (known) existence of such algebras over the fields
$F_{C_i}$, modify the construction to
produce a class in $\H^2(\O_{C,\S})\oplus\Gamma$, and then apply
Theorem \ref{map} and Theorem \ref{preservesindex} to prove existence over $F$.

By \cite[Theorem 4.7]{BMT}, if $F$ admits a smooth $X$
then there exist noncrossed product division algebras
over $F_C$ of period and index as low as
$\ell^{s+1}$ if $r=0$, and $\ell^{2r+1}$ if $r>0$.
The resulting Brauer class has the form 
$\alpha_{C}+(\pi)\cdot\theta_{C}\in\H^2(F_{C})$,
where $\alpha_{C}\in\H^2(\kappa(C))$ and $\theta_{C}\in\H^1(\kappa(C),-1)$.
A look at the construction, which proceeds exactly as in \cite[Theorem 1]{Br95},
shows we may pre-assign values at any finite set of points of $C$. 
Thus we may thus produce a noncrossed product $F_{C_1}$-division
algebra with class $\gamma_{C_1}=\alpha_{C_1}+(\pi)\cdot\theta_{C_1}$ 
of the desired period and index, and elements $\gamma_{C_i}=\alpha_{C_i}+(\pi)\cdot\theta_{C_i}$ for $i>1$ 
that glue as in Lemma~\ref{gluing}.  Let
$$
\gamma_C=
(\alpha_{C_1}+(\pi)\cdot\theta_{C_1},\dots,
\alpha_{C_m}+(\pi)\cdot\theta_{C_m})\;\in\H^2(\O_{C,\S})\oplus(\pi)\cdot\Gamma\leq\bigoplus_i\H^2(F_{C_i})
$$
Then $\gamma_C$ lifts to $\gamma=\lambda(\gamma_C)\in\H^2(F)$ by Theorem \ref{map}. 
By Theorem \ref{preservesindex}, $\ind(\gamma)=\ind(\gamma_C)=\ell^{s+1}$ if $r=0$, and $\ell^{2r+1}$ if $r\neq 0$.
It is clear that the $F$-division algebra $D$ associated to $\gamma$ is a noncrossed
product, since any Galois maximal subfield of $D$ over $F$ could be constructed over $F_{C_1}$,
contradicting the fact that $D\otimes_F F_{C_1}$ is a noncrossed product $F_{C_1}$-division algebra.

\hfill $\blacksquare$

\Theorem
Let $F/\Q_p$ be a finitely generated field extension of transcendence degree one,
and let $\ell\neq p$ be a prime.  Then there exist indecomposable $F$-division algebras
of (period,index)$=(\ell^a,\ell^b)$, for any numbers $a$ and $b$ satisfying
$1\leq a\leq b\leq 2a-1$.
\rm

\Pf
Let $X$, $C$, $C_i$, and $\S$ be as in Theorem \ref{noncrossed}.
The construction over $F_{C_i}$ is exactly as in \cite[Proposition 4.2]{BMT}
and \cite{Br96b},
and we merely have to observe that we may assume all of the data
in the constructed class $\gamma_{C_i}=\alpha_{C_i}+(\pi)\cdot\theta_{C_i}\in\H^2(F_{C_i})$
is trivial at
the singular points $\S\cap C_i$, so that by Lemma~\ref{gluing}, we may construct
a class $\gamma_C=\alpha_C+(\pi)\cdot\theta_C$ in 
$\H^2(\O_{C,\S})\oplus(\pi)\cdot\Gamma
\leq\bigoplus_i\H^2(F_{C_i})$ whose $i$-th component is $\alpha_{C_i}+(\pi)\cdot\theta_{C_i}$.
This class lifts to a class $\gamma=\lambda(\gamma_C)$ by Theorem \ref{map}, 
and $\ind(\gamma)=\ind(\gamma_C)$ by Theorem \ref{preservesindex}.
Since the indexes are the same, the division algebra $D$ associated to $\gamma$
is indecomposable, since any decomposition would extend to
$D_{C_i}=D\otimes_F F_{C_i}$, contradicting the construction of $\gamma_{C_i}$.

\hfill $\blacksquare$

\bibliographystyle{abbrv} 
\bibliography{hnx.bib}

\end{document}